\documentclass{article}
\usepackage{amsmath,amsfonts,amsthm}
\newtheorem{thm}{Theorem}[section]
\newtheorem{pro}[thm]{Proposition}
\newtheorem{cor}[thm]{Corollary}
\newtheorem{lem}[thm]{Lemma}
\newtheorem{conj}[thm]{Conjecture}

\newtheorem{claim}[thm]{Claim}
\newtheorem{ex}[thm]{Example}
\newtheorem{rem}[thm]{Remark}
\newtheorem{defn}[thm]{Definition}
\newtheorem{quot}[]{Result}
\newcommand{\dime}{\operatorname{dim}}

\def\div{\raise 1pt \hbox{\big|}}

\begin{document}

\title{On Commuting and Non-Commuting Complexes}
\author{Jonathan Pakianathan and Erg\"un Yal\c c\i n}
\maketitle

\begin{abstract}
In this paper we study various simplicial complexes associated 
to the commutative structure of a finite group $G$. 
We define $NC(G)$ (resp. $C(G)$) as the complex associated 
to the poset of pairwise non-commuting (resp. commuting)
sets of nontrivial elements in $G$. 

We observe that $NC(G)$ has only 
one positive dimensional connected component, which we call $BNC(G)$,
and we prove that $BNC(G)$ is simply connected.

Our main result is a simplicial 
decomposition formula for $BNC(G)$ which follows
from a result of A. Bj\"orner, M. Wachs and V. Welker on
inflated simplicial complexes (See \cite{Bj}.)
As a corollary we obtain that if $G$ has a nontrivial center or if $G$ has 
odd order, then the homology group $H_{n-1}(BNC(G))$ is nontrivial for every
$n$ such that $G$ has a maximal noncommuting set of order $n$.

We discuss the duality between $NC(G)$ and $C(G)$, and between
their $p$-local versions $NC_p(G)$ and $C_p(G)$. We observe
that $C_p(G)$ is homotopy equivalent to the Quillen complexes $A_p(G)$,
and obtain some interesting results for $NC_p(G)$ using this duality.

Finally, we study the family of groups where the commutative
relation is transitive, and show that in this case,  
$BNC(G)$ is shellable. As a consequence we derive
some group theoretical formulas for the orders of maximal
non-commuting sets.

\end{abstract}

\section{Introduction}

Given a finite group $G$, one defines a non-commuting set 
to be a set of elements $\{g_1, \dots, g_n\}$ such that 
$g_i$ does not commute with $g_j$ for $i \neq j$. The sizes of maximal
non-commuting sets in a group are interesting invariants of the group.
In particular, the largest integer $n$ such that the group $G$ has a non-commuting
set of order $n$, which is denoted by $nc(G)$, is known to be closely
related to other invariants of $G$. For example if $k(G)$ is the size of 
the largest conjugacy class in $G$ then 
$$
k(G) \leq 4(nc(G))^2
\quad \quad \text{(See \cite{P}.)} 
$$
Also if we define $cc(G)$ to be the minimal 
number of abelian subgroups of $G$ that covers $G$, 
then I. M. Isaacs (see \cite{J}) has shown that  
$$
nc(G) \leq cc(G) \leq (nc(G)!)^2.
$$

Confirming a conjecture of P. Erd\"os, 
L. Pyber \cite{P} has also shown that there is a positive constant $c$ 
such that 
$$
cc(G) \leq |G:Z(G)| \leq c^{nc(G)}
$$
for all groups $G$. Another interesting place where the invariant
$nc(G)$ appears is in the computation of the cohomology length of
extra-special $p$-groups (See \cite{Y}.) 

In this paper we study the topology of certain 
complexes associated to the poset of non-commuting sets in a group $G$. 
Let $NC(G)$ be the complex whose vertices are just the nontrivial 
elements of the group $G$ and whose faces are the noncommuting sets 
in $G$. The central elements form point components in this complex 
and are not so interesting. So, we look at the subcomplex $BNC(G)$
formed by non-central elements of $G$. We show

\begin{quot}[~\ref{thm: sc}]
If $G$ is a nonabelian group, then $BNC(G)$ is simply-connected. 
\end{quot}

In general we also notice that $BNC(G)$ is equipped with a free 
$Z(G)$-action where $Z(G)$ is the center of $G$. It is also 
equipped in general with a $\mathbb{Z}/2\mathbb{Z}$-action 
whose fixed point set is exactly $BNC_2(G)$, the corresponding complex 
where we use only the elements of order 2 (involutions). Thus if $G$ is an odd 
order group or if $G$ has nontrivial center, then the Euler characteristic of 
$BNC(G)$ is not $1$ and so it is not contractible. 

On a more refined level, we use a recent simplicial decomposition 
result of A. Bj\"orner, M. Wachs and V. Welker to show that  
there is a simplicial complex $S$, called the core of $BNC(G)$, 
so that the following decomposition formula holds:

\begin{quot}[~\ref{thm: core}]
If $G$ is a finite nonabelian group and $S$ is the core of $BNC(G)$, then 
$$
BNC(G) \simeq S \vee \bigvee_{F \in S} 
[Susp^{|F|}(Lk F)]^{\gamma(F)}
$$
where the $F$ are the faces of $S$, $Lk$ stands for the link of a face, 
$Susp^k$ stands for a $k$-fold suspension and $\vee$ stands for wedge 
product. The number $\gamma(F)=\prod_{[x] \subseteq F}(m_x - 1)$ 
where $m_x$ is the size of the centralizer class $[x]$.
\end{quot}

It is clear from this decomposition formula that when the core of $BNC(G)$ 
is contractible, then  $BNC(G)$ is a wedge of suspensions of spaces 
and hence has a trivial ring structure on its cohomology. This is
true, for example, when $G=\Sigma_p$, the symmetric group on 
$p$ letters, for some prime $p$. 

The following is an important consequence of the above decomposition
formula:

\begin{quot}[~\ref{cor: dimension}]
\label{result:dimension}
Let $G$ be a finite nonabelian group, and let $S_s$ 
denote the set of maximal non-commuting sets in 
$G$ of size $s$. Then, for $s >1$,
$$
rk(H_{s-1}(BNC(G))) \geq \sum_{F \in S_s} [\prod_{x \in F}
(1-\frac{1}{m_x})] 
$$
where $m_x$ is the size of the centralizer class containing $x$.

In particular, if $G$ is an odd order group 
or if $G$ has a nontrivial center $(|G|\neq 2)$, then 
$$ 
{\tilde H}_{s-1}(NC(G)) \neq 0 
$$ 
whenever $G$ has a maximal non-commuting set of size $s$.
\end{quot}

There is also a $p$-local version of this theorem which, in particular,
gives that ${\tilde H}_{s-1}(NC_p(G)) \neq 0$ whenever $G$ has a 
maximal non-commuting $p$-set of size $s$ and $p$ is an odd prime.
For $p=2$, the same is true under the condition $2| |Z(G)|$
and $|G| \neq 2$.
Observe that this result has a striking formal similarity (in terms
of their conclusions) to the following theorem proved by Quillen (theorem
12.1 in \cite{Q}):

\begin{thm}[Quillen] If $G$ is a finite solvable group having no nontrivial 
normal $p$-subgroup, then 
$${\tilde H}_{s-1} (A_p(G)) \neq 0 $$ 
whenever $G$ has a maximal elementary abelian $p$-group of rank $s$.
\end{thm}

A nice consequence of Result~\ref{result:dimension} can be stated as follows:
If $G$ is a group of odd order or a group with nontrivial
center such that $BNC(G)$ is spherical, i.e., homotopy
equivalent to a wedge of equal dimensional
spheres, then all maximal non-commuting sets in $G$ have the same size.

A natural question to ask is: For which groups is $BNC(G)$ spherical?
As a partial answer, we show that if $G$
is a group  where the commutation relation is transitive, then $BNC(G)$
is spherical. We give examples of such groups and compute $BNC(G)$ for
these groups.

Observe that one could also define a commuting complex $C(G)$ analogous to 
the way we defined $NC(G)$ by making the faces consist of commuting sets 
of elements instead of non-commuting sets. However, this definition
does not provide us with new complexes. For example, $C_p(G)$, the commuting 
complex formed by the elements of prime order $p$, is easily shown to 
be $G$-homotopy equivalent to Quillen's complex $A_p(G)$. However, 
the definition
of $C_p(G)$ helps us to see a duality between $NC_p(G)$ and $C_p(G)$
where $NC_p(G)$ is the subcomplex of $NC(G)$ spanned by the vertices 
which correspond to elements of order $p$. Using a result 
of Quillen on $A_p(G)$, we obtain:

\begin{quot}[~\ref{thm: Quillen}]
Let $G$ be a group and $p$ a prime with $p | |G|$. Pick a Sylow $p$-subgroup 
$P$ of $G$ and define $N$ to be the subgroup generated by 
the normalizers $N_G(H)$ where $H$ runs over all the nontrivial subgroups of 
$P$. 

Then $NC_p(G)$ is $(|G:N|-2)$-connected. In fact $NC_p(G)$ is the 
$|G:N|$-fold join of a complex $S$ with itself where $S$ 
is ``dual'' to a path-component of $A_p(G)$. 
\end{quot}  

Finally under suitable conditions, $BNC(G)$ is shellable and this yields  
some combinatorial identities. As an application we obtain

\begin{quot}[~\ref{cor: Shell}]
Let $G$ be a nonabelian group with a transitive commuting relation,
i.e., if $[g,h]=[h,k]=1$, then $[g,h]=1$ for every noncentral $g,h,k \in G$.
Then, 
$$
nc(G)(nc(G)-1) + |G|(|G|-m) - 2(nc(G)-1)(|G|-|Z(G)|) \geq 0.
$$
where $m$ denotes the number of conjugacy classes in $G$.
\end{quot}

\section{Background}

We start the section with a discussion of complexes associated
with posets of subgroups of a group $G$. For a complete account of
these well known results, see chapter 6 in \cite{Ben}.

Given a finite poset $(P, \leq)$, one can construct a 
simplicial complex $|P|$ out of it by defining the $n$-simplices of $|P|$ 
to be chains in $P$ of the form $p_0 < p_1 < \dots < p_n$. 
This is called the simplicial
realization of the poset $P$.

Furthermore any map of posets $f: (P_1, \leq_1) \rightarrow (P_2, \leq_2)$ 
(map of posets means $x_1 \leq_1 x_2 \Rightarrow f(x_1) \leq_2 f(x_2)$)
yields a simplicial map between $|P_1|$ and $|P_2|$  
and hence one has in general a covariant functor from the category 
of finite posets to the category of finite simplicial complexes and 
simplicial maps.
Thus if a (finite) group $G$ acts on a poset $P$ via poset maps (we say 
$P$ is a $G$-poset) then $G$ will act on $|P|$ simplicially.

Brown, Quillen, Webb, Bouc, Th\' evenaz and many others 
constructed many finite $G$-simplicial complexes associated to a group 
$G$ and used them to study the group $G$ and its cohomology.
In particular, the following posets of subgroups of $G$ have been
studied extensively: \\

\noindent
(a) The poset $s_p(G)$ of nontrivial $p$-subgroups of $G$. \\ 
\noindent
(b) The poset $a_p(G)$ of nontrivial elementary abelian $p$-subgroups of 
$G$. \\
\noindent
(c) The poset $b_p(G)$ of nontrivial $p$-radical subgroups of $G$. 
(Recall a $p$-radical subgroup of $G$ is a $p$-subgroup $P$ of 
$G$ such that $PN_G(P)/P$ has no nontrivial normal $p$-subgroups.) \\

$G$ acts on each of these posets by conjugation and thus from each 
of these $G$-posets, one gets a $G$-simplicial complex 
$S_p(G)$, $A_p(G)$ and $B_p(G)$ respectively. 
$S_p(G)$ is usually called the Brown complex of $G$ and 
$A_p(G)$ is usually called the Quillen complex of $G$ where the dependence 
on the prime $p$ is understood. Notice again that the trivial subgroup is 
not included in any of these posets, since if it were the resulting complex 
would be a cone and hence trivially contractible.

It was shown via work of Quillen and Th\' evenaz, that $S_p(G)$ and 
$A_p(G)$ are $G$-homotopy equivalent and via work of Bouc and Th\' evenaz 
that $B_p(G)$ and $S_p(G)$ are $G$-homotopy equivalent. 
Thus, in a sense, these three $G$-complexes capture the same information. 

Recall the following elementary yet very important lemma 
(see \cite{Ben}):

\begin{lem}
\label{lem: basic}
If $f_0, f_1 : P_1 \rightarrow P_2$ are two maps of posets such that 
$f_0(x) \leq f_1(x)$ for all $x \in P_1$ then the simplicial maps 
induced by $f_0$ and $f_1$ from $|P_1|$ to $|P_2|$ are 
homotopic.
\end{lem}

Using this, Quillen made the following observation, if $P_0$ is a nontrivial  
normal $p$-subgroup of $G$, then we may define a poset map 
$f: s_p(G) \rightarrow s_p(G)$ by $f(P)=P_0P$ and by the lemma above, 
$f$ would be homotopic to the identity map, but on the other hand since 
$f(P)$ contains $P_0$ always, again by the lemma, 
$f$ is also homotopic to a constant map. Thus we see that $S_p(G)$ is 
contractible in this case. Quillen then conjectured

\begin{conj}[Quillen]
If $G$ is a finite group, $S_p(G)$ is contractible if and only if 
$G$ has a nontrivial normal $p$-subgroup.
\end{conj}

He proved his conjecture in the case that $G$ is solvable but the 
general conjecture remains open. Notice though that 
if $S_p(G)$ is $G$-homotopy equivalent to a point space then this does imply  
$G$ contains a nontrivial normal $p$-subgroup since in this case $S_p(G)^G$ 
is homotopy equivalent 
to a point which means in particular $S_p(G)^G$ is not empty, 
yielding a nontrivial normal $p$-subgroup.

The purpose of this paper is to introduce some simplicial 
complexes associated to 
elements of a group rather than subgroups of a group and use these 
to give a different perspective on some of the complexes above.

For this purpose we make the following definitions:

\begin{defn}
Let $G$ be a group. 
Define a simplicial complex $C(G)$ by declaring a $n$-simplex in this 
complex to be a collection $[g_0,g_1,\dots,g_n]$ of distinct nontrivial 
elements of $G$, which pairwise commute.

Similarly define a simplicial complex $NC(G)$ by declaring a 
$n$-simplex to be a collection $[g_0,g_1,\dots,g_n]$ of nontrivial 
elements of $G$, which pairwise do not commute.
\end{defn}

It is trivial to verify that the above definition, does indeed define 
complexes on which $G$ acts simplicially by conjugation. 

Usually when one studies simplicial group actions, it is nice to have 
admissible actions, i.e., actions where if an element of $G$ fixes a 
simplex, it actually fixes it pointwise. Although $C(G)$ and $NC(G)$ 
are not admissible in general, one can easily fix this by taking 
a barycentric subdivision. The resulting complex is of course $G$-homotopy 
equivalent to the original, however it now is the realization of a poset.

Thus if we let $PC(G)$ be the barycentric subdivision of $C(G)$, it 
corresponds to the realization of the poset consisting of 
subsets of nontrivial, pairwise commuting elements 
of $G$, ordered by inclusion. 
Similarly if we let $PNC(G)$ be the barycentric subdivision of $NC(G)$, 
it corresponds to the realization of the poset consisting of 
subsets of nontrivial, pairwise non-commuting elements of $G$, 
ordered by inclusion. 

Depending on the situation, one uses either the barycentric subdivision 
or the original. For purposes of understanding the topology, the original 
is easier but for studying the $G$-action, the subdivision is easier.

Of course, we will want to work a prime at a time also, so we introduce 
the following $p$-local versions of $C(G)$ and $NC(G)$.

\begin{defn}
Given a group $G$ and a prime $p$, let
$C_p(G)$ be the subcomplex of $C(G)$ where the simplices consist 
of sets of nontrivial, pairwise commuting elements of order $p$.

Similarly let $NC_p(G)$ be the subcomplex of $NC(G)$ where  
the simplices consist 
of sets of nontrivial, pairwise non-commuting elements of order $p$  
\end{defn}

Of course the same comments about the $G$-action and the barycentric 
subdivision apply to these $p$-local versions. 

Our first order of business is to see that the commuting complexes 
$C(G)$ and $C_p(G)$ are nothing new. We will find the following 
standard lemma useful for this purpose (see \cite{Ben}):

\begin{lem}
\label{lem: Ghom}
If $f$ is a $G$-map between admissible 
$G$-simplicial complexes $X$ and $Y$ with 
the property that for all subgroups $H \leq G$, f restricts to 
an ordinary homotopy equivalence between $X^H$ and $Y^H$ (recall 
$X^H$ is the subcomplex of $X$ which consists of elements fixed 
pointwise by $H$), then 
$f$ is a $G$-homotopy equivalence, i.e., there is a $G$-map 
$g: Y \rightarrow X$ such that $f \circ g$ and $g \circ f$ are 
$G$-homotopic to identity maps. 
\end{lem}

\begin{thm}
\label{thm: OldCNewC}
Let $G$ be a finite group, 
then $C(G)$ is $G$-homotopy equivalent to the simplicial realization 
of the poset $A(G)$ 
of nontrivial abelian subgroups of $G$, ordered by inclusion 
and acted on by conjugation. 

Furthermore if $p$ is a prime, then 
$C_p(G)$ is $G$-homotopy equivalent to $A_p(G)$ (and thus $S_p(G)$ and 
$B_p(G)$.)
\end{thm}
\begin{proof}
First we will show homotopy equivalence and remark on 
$G$-homotopy equivalence later. 

We work with $PC(G)$, the barycentric subdivision. Notice that the 
associated poset of $PC(G)$ contains the poset $A(G)$ of nontrivial abelian 
subgroups of $G$ as a subposet, they are merely the commuting sets 
whose elements actually form an abelian subgroup (minus identity). 
Let $i: A(G) \rightarrow PC(G)$ denote this inclusion. 

We now define a poset map $r: PC(G) \rightarrow A(G)$ as follows: 
If $S$ is a set of nontrivial, pairwise commuting elements of $G$, 
then $<S>$, the subgroup generated by $S$ will be a nontrivial abelian 
subgroup of $G$, thus we can set $r(S)=\{<S>-1\}$. It is obvious  
that $r$ is indeed a poset map, and that $S \subset r(S)$ and so 
$i \circ r$ is homotopic to the identity map of $PC(G)$ by 
lemma~\ref{lem: basic}. Furthermore it is clear that $r \circ i = Id$ 
and so $r$ is a deformation retraction of $PC(G)$ onto $A(G)$. 

Thus $PC(G)$ is homotopy equivalent to $A(G)$.  
To see this is a $G$-homotopy equivalence, we just need to note 
that $r$ is indeed a $G$-map, and maps a commuting set invariant under 
conjugation by 
a subgroup $H$ into a subgroup invariant under conjugation 
by $H$ and thus induces 
a homotopy equivalence between $PC(G)^H$ and $A(G)^H$ for 
any subgroup $H$. Thus $r$ is indeed a $G$-homotopy equivalence 
by lemma~\ref{lem: Ghom}.

The $p$-local version follows exactly in the same manner, once one 
notes that the subgroup generated by a commuting set of elements of order $p$ 
is an elementary abelian $p$-subgroup.  
\end{proof}

Thus we see from theorem~\ref{thm: OldCNewC}, that the commuting complexes 
at a prime $p$ are basically the $A_p(G)$ in disguise. 
However for the rest of the paper, we look at the non-commuting complexes 
and we will see that they are quite different, and in some sense dual 
to the commuting ones.  However before doing that, we conclude this 
section by looking at a few more properties of the commuting complex.

Recall the following important proposition of Quillen \cite{Q}:

\begin{pro}
\label{pro: basic}
If $f: X \rightarrow Y$ is a map of posets, and $y \in Y$ we define
\begin{align*}
\begin{split}
f/y &= \{ x \in X | f(x) \leq y \} \\
y\backslash f &= \{ x \in X | f(x) \geq y \}.
\end{split}
\end{align*}
Then if $f/y$ is contractible for all $y \in Y$ (respectively 
$y\backslash f$ is contractible for all $y \in Y$) then $f$ is a homotopy 
equivalence between $|X|$ and $|Y|$.
\end{pro}

Using this we will prove:

\begin{pro}
If $G$ is a group 
with nontrivial center then $A(G)$ and hence $C(G)$ is contractible.
Moreover, $A(G)$ is homotopy equivalent to $Nil(G)$ where $Nil(G)$ is the poset 
of nontrivial nilpotent subgroups of $G$. 
\end{pro}
\begin{proof}
If $G$ has a nontrivial center $Z(G)$, then $A(G)$ is conically 
contractible via $A \leq AZ(G) \geq Z(G)$ for any abelian subgroup 
$A$ of $G$. Thus $C(G)$ is also contractible as it is homotopy equivalent 
to $A(G)$.

Let $i: A(G) \rightarrow Nil(G)$ be the natural inclusion of posets.
Take $N \in Nil(G)$ and let us look at $i/N = \{B \in A(G) | B \subseteq N \} 
=A(N)$. However since $N$ is nilpotent, it has a nontrivial center $Z$ 
and hence $A(N)$ is conically contractible. 
Thus by proposition~\ref{pro: basic} the result follows.

\end{proof}

It is natural to ask if:

\begin{conj}
$C(G)$ is contractible if and only if $G$ has a nontrivial center.
\end{conj}

\section{Non-commuting complexes}

Fix a group $G$, let us look at $NC(G)$. The first thing we notice 
is that any nontrivial central element in $G$ gives us a point component 
in $NC(G)$ and hence is not so interesting. Thus we define:

\begin{defn}
$BNC(G)$ is the subcomplex of $NC(G)$ consisting of those simplices 
of $NC(G)$ which are made out of noncentral elements.
Thus $BNC(G)$ is empty if $G$ is an abelian group.
\end{defn}

The first thing we will show is that if $G$ is a nonabelian group, 
(so that $BNC(G)$ is nonempty) then $BNC(G)$ is not only path-connected 
but it is simply-connected. Notice also that this means the general 
picture of $NC(G)$ is as a union of components, with at most one 
component of positive dimension and this is $BNC(G)$ and it 
is simply-connected. Also notice that $BNC(G)$ is invariant under the 
conjugation $G$-action, and the point components of $NC(G)$ are just fixed 
by the $G$-action as they correspond to central elements.

\begin{thm}
\label{thm: sc}
If $G$ is a nonabelian group, then $BNC(G)$ is a simply-connected 
$G$-simplicial complex. 
\end{thm}
\begin{proof}
First we show that $BNC(G)$ is path-connected. 
Take any two vertices in $BNC(G)$, call them $g_0$ and $g_1$, then 
these are two noncentral elements of $G$, thus their centralizer 
groups $C(g_0)$ and $C(g_1)$ are proper subgroups of $G$. 

It is easy to check that no group is 
the union of two proper subgroups for suppose $G=H \cup K$ 
where $H, K$ are proper subgroups of $G$. Then we can find $h \in G-K$ 
(it follows $h \in H$) and $k \in G-H$ (hence $k \in K$). 
Then $hk$ is not in $H$ as $h \in H$ and $k \notin H$. Similarly 
$hk \notin K$ so $hk \notin H \cup K = G$ which is an obvious contradiction. 
Thus no group is the union of two proper subgroups. 

Thus we conclude that $C(g_0) \cup C(g_1) \neq G$ and so we can find 
an element $w$ which does not commute with either $g_0$ or $g_1$ and 
so the vertices $g_0$ and $g_1$ are joined by an edge path $[g_0, w] + 
[w, g_1]$. (The $+$ stands for concatenation.) 
Thus we see $BNC(G)$ is path-connected. In fact, 
any two vertices of $BNC(G)$ can be connected by an edge path 
involving at most two edges of $BNC(G)$.

To show it is simply-connected, we argue by contradiction. If it was not 
simply connected, 
then there would be a simple edge loop which did not contract, 
i.e., did not bound a suitable union of 2-simplices. 
(a simple edge loop is formed by edges of the simplex 
and is of the form $L=[e_0, e_1] + [e_1, e_2] + \dots + [e_{n-1},e_n]$ where 
all the $e_i$ are distinct except $e_0=e_n$.) 

Take such a loop $L$ with minimal size $n$.  (Notice $n$ is just the number 
of edges involved in the loop.)

Since we are in a simplicial complex, certainly $n \geq 3$. 

Suppose $n>5$, then $e_3$ can be connected to $e_0$ by an edge path $E$  
involving at most two edges by our previous comments. This edge path $E$ 
breaks our simple edge loop into two edge loops of smaller size which hence 
must contract since our loop was minimal. However, then it is clear 
that our loop contracts which is a contradiction so $n \leq 5$. 

So we see $3 \leq n \leq 5$ so we have three cases to consider:

\noindent
{\bf (a) $n=3$:} \\
Here $L=[e_0,e_1] + [e_1,e_2] + [e_2,e_3]$ with $e_3=e_0$.
But then it is easy to see $\{e_0, e_1, e_2\}$ is a set of pairwise 
non-commuting elements and so gives us a 2-simplex $[e_0,e_1,e_2]$ 
in $BNC(G)$ which bounds the loop which gives a contradiction.

\noindent
{\bf (b) $n=4$:} \\
Here $L=[e_0,e_1] + [e_1,e_2] + [e_2,e_3] + [e_3,e_4]$ with 
$e_4=e_0$. Thus $L$ forms a square. Notice by the simplicity of $L$, 
the diagonally opposite vertices in the square must not be joined 
by an edge in $BNC(G)$, i.e., they must commute, thus 
$e_0$ commutes with $e_2$ and $e_1$ commutes with $e_3$.

Since $e_0$ and $e_1$ do not commute, $\{e_0, e_1, e_0e_1 \}$ is 
a set of mutually non-commuting elements and so forms a 2-simplex 
of $BNC(G)$. Since $e_2$ commutes with $e_0$ but not with $e_1$, 
it does not commute with $e_0e_1$ and thus $\{e_0e_1, e_1, e_2\}$ 
also is a 2-simplex in $BNC(G)$. Similar arguments show that 
$\{e_0e_1, e_0, e_3 \}$ and $\{e_0e_1, e_2, e_3 \}$ form 
2-simplices in $BNC(G)$. The union of the four $2$-simplices mentioned 
in this paragraph, bound our loop giving us our contradiction. 
Thus we are reduced to the final case:

\noindent
{\bf (c) $n=5$:} \\
Here $L=[e_0,e_1] + [e_1,e_2] + [e_2,e_3] + [e_3,e_4] + [e_4,e_5]$ 
with $e_5=e_0$. Thus $L$ forms a pentagon, and again by simplicity 
of $L$, nonadjacent vertices in the pentagon cannot be joined by an 
edge in $BNC(G)$, thus they must commute. Similar arguments as for 
the $n=4$ case yield that $[e_0, e_1, e_0e_1]$, $[e_0e_1, e_1, e_2]$ 
and $[e_0e_1, e_0, e_4]$ are $2$-simplices in $BNC(G)$ and the union of 
these three simplices contract our loop $L$ into one of length four 
namely $[e_0e_1, e_2] + [e_2,e_3] + [e_3,e_4] + [e_4, e_0e_1]$ which 
by our previous cases, must contract thus yielding the final 
contradiction. 
\end{proof}

From theorem~\ref{thm: sc}, we see that $BNC(G)$ is simply-connected for 
any nonabelian group $G$. One might ask if it is contractible? 
The answer is no in general although there are groups $G$ where it is 
contractible. We look at these things next:

\begin{pro}
\label{pro: oddorder}
For a general finite nonabelian group $G$, 
the center $Z(G)$ of $G$ acts freely on $BNC(G)$ by left multiplication 
and hence $|Z(G)|$ divides the Euler characteristic of $BNC(G)$. 

For a group of odd order, the simplicial map $A$ which  
maps a vertex $g$ to $g^{-1}$ is a fixed point free map on 
$BNC(G)$ and on $NC(G)$ of order $2$.
Thus the Euler characteristic of both $BNC(G)$ and $NC(G)$ 
is even in this case. 

Thus if $BNC(G)$ is $\mathbb{F}$-acyclic for some field $\mathbb{F}$, then 
$G$ must have trivial center and be of even order. 
\end{pro}
\begin{proof}
The remarks about Euler characteristics follow from the fact 
that if a finite group 
$H$ acts freely on a space where the Euler characteristic is 
defined, then $|H|$ must divide 
the Euler characteristic. So we will concentrate mainly on 
finding such actions.

First for the action of $Z(G)$. $a \in Z(G)$ acts by taking a 
simplex $[g_0,\dots,g_n]$ to a simplex $[ag_0,\dots,ag_n]$. 
Notice this is well-defined since $ag_i$ is noncentral if $g_i$ is   
noncentral and since $ag_i$ commutes with $ag_j$ if and only if 
$g_i$ commutes with $g_j$.

Furthermore if $a$ is not the identity element, this does not fix 
any simplex $[g_0, \dots, g_n]$ since if the set 
$\{g_0,\dots,g_n\}$ equals the set $\{ag_0,\dots,ag_n\}$ 
then $ag_0$ is one of the $g_j$'s but $ag_0$ commutes with $g_0$ so 
it would have to be $g_0$ but $ag_0=g_0$ gives $a=1$, a contradiction.

Thus this action of nonidentity central $a$ does not fix any simplex 
and so we get a free action of $Z(G)$ on $BNC(G)$. 

Now for the action of $A$, first note that $A$ is well-defined since 
if $[g_0,\dots,g_n]$ is a set of mutually non-commuting elements of $G$, 
so is $[g_0^{-1}, \dots ,g_n^{-1}]$. Clearly 
$A \circ A = Id$. Furthermore, if the two sets above are equal, 
then $g_0^{-1}$ would have to be one 
of the $g_j$. But since $g_0^{-1}$ commutes with $g_0$ it would have to 
be $g_0$, i.e., $g_0$ would have to have order $2$. Similarly, 
all the $g_i$ would have to have order 2. 
Thus in a group of odd order, $A$ would not fix any simplex of $BNC(G)$ 
or $NC(G)$, and hence would not have any fixed points.

The final comment is to recall that if $BNC(G)$ were $\mathbb{F}$-acyclic, 
its Euler characteristic would be 1 and hence the center of $G$ would 
be trivial and $G$ would have to have even order by the facts we 
have shown above.
\end{proof}

Proposition~\ref{pro: oddorder} shows that $BNC(G)$ is not contractible 
if $G$ is of odd order or if $G$ has a nontrivial center, for example 
if $G$ were nilpotent. There is a corresponding $p$-local version 
which we state next:

\begin{pro}
If $2 | |Z(G)|$ or if $G$ has odd order then $NC_2(G)$ is not contractible, 
in fact the Euler characteristic is even.
$NC_p(G)$ is never contractible for any odd prime $p$, in fact it always 
has even Euler characteristic.
\end{pro}
\begin{proof}
Follows from the proof of proposition~\ref{pro: oddorder}, once 
we note that left multiplication by a central element of order $2$
takes the subcomplex $NC_2(G)$ of $NC(G)$ to itself and that the 
map $A$ maps $NC_p(G)$ into itself, as the inverse of an element 
has the same order as the element. For odd primes $p$, $A$ is fixed 
point free on $NC_p(G)$ always as no elements of order 2 are involved 
in $NC_p(G)$. 
\end{proof} 

The fact that 
$BNC(G)$ can be contractible sometimes is seen in the next proposition.

\begin{pro}
If $G$ is a nonabelian group with a self-centralizing involution i.e., an 
element $x$ of order 2 such that $C(x)=\{1,x\}$, then 
$BNC(G)$ is contractible. In fact, $BNC(G)=NC(G)$ in this case. 

Thus for example, $NC(\Sigma_3)=BNC(\Sigma_3)$ is contractible 
where $\Sigma_3$ is the symmetric group on 3 letters.
\end{pro}
\begin{proof}
Since $x$ does not commute with any nontrivial element, the center of 
$G$ is trivial and $NC(G)=BNC(G)$. Furthermore, it is clear that 
$BNC(G)$ is a cone with $x$ as vertex and hence is contractible.
\end{proof}

To help show $NC(G)$ of a group is not contractible, we note the following 
observation which uses Smith theory. (See \cite{Ben}).

\begin{pro}
\label{pro: 2local}
If $G$ is a group and $NC(G)$ is $\mathbb{F}_2$-acyclic where 
$\mathbb{F}_2$ is the field with two elements, then 
$NC_2(G)$ is also $\mathbb{F}_2$-acyclic. Furthermore one always has 
$\chi(NC(G))=\chi(NC_2(G))$ mod $2$.
\end{pro}
\begin{proof}
We first recall that the map $A$ from proposition~\ref{pro: oddorder} 
has order 2 as a map of $NC(G)$. However 
it might have fixed points, in fact from the proof of 
that proposition, we see that $A$ fixes a simplex $[g_0,\dots,g_n]$ 
of $NC(G)$ if and only if each element $g_i$ has order $2$ and 
it fixes the simplex pointwise. Thus the fixed point set of $A$ 
on $NC(G)$ is nothing other than $NC_2(G)$ the $2$-local non-commuting 
complex for $G$. Since $A$ has order 2, we can apply Smith Theory 
to finish the proof of the first statement of the proposition.
The identity on the Euler characteristics follows once we note 
that under the action of $A$, 
the cells of $NC(G)$ break up into free orbits and cells which 
are fixed by $A$ and the fixed cells exactly form $NC_2(G)$. 
\end{proof} 

We observed earlier
that if $BNC(G)$ is contractible then the center of $G$ is trivial, 
i.e. $BNC(G)=NC(G)$, hence, by proposition~\ref{pro: 2local}, 
$NC_2(G)$ is $\mathbb{F}_2$-acyclic (and in particular nonempty).

\begin{defn}
Let $G$ be a finite group. We define $nc(G)$ to be the 
maximum size of a set of pairwise non-commuting elements in $G$. 
Thus $nc(G) - 1$ is the dimension of $NC(G)$.
\end{defn}

We now compute a general class of examples, the Frobenius groups. 
Recall a group $G$ is a Frobenius group if it has a proper nontrivial 
subgroup $H$ with the property that $H \cap H^g = 1$ if 
$g \in G-H$. $H$ is called the Frobenius complement of $G$. 
Frobenius showed the existence of a normal subgroup $K$ such 
that $K=G-\cup_{g \in G}(H^g-\{1\})$. Thus $G$ is a split
extension of $K$ by $H$, i.e., $G=K \times_{\phi} H$, 
for some homomorphism $\phi: H \rightarrow Aut(K)$.
This $K$ is called the Frobenius kernel of $G$.

We have:

\begin{pro}
\label{pro: Frob}
If $G$ is a Frobenius group with Frobenius kernel $K$ and 
Frobenius complement $H$, then 
$$
NC(G)=NC(K) * NC(H)^{|K|}
$$
where $*$ stands for simplicial join and the superscript $|K|$ means that 
$NC(K)$ is joined repeatedly with $|K|$ many copies of $NC(H)$.

It also follows that $nc(G)=nc(K) + |K|nc(H)$.

Finally if both $H$ and $K$ are abelian, then $NC(G)$ is 
homotopy equivalent to a wedge of $(|K|-2)(|H|-2)^{|K|}$ many $|K|$-spheres.
\end{pro}
\begin{proof}
First, from the condition that $H \cap H^g = 1$ for 
$g \in G-H$ we see that no nonidentity element of $H$ commutes 
with anything outside of $H$. Thus conjugating the picture, 
no nonidentity element of any $H^g$ commutes with anything 
outside $H^g$. Thus if $H^{g_1}, H^{g_2},\dots,H^{g_m}$ is a 
complete list of the conjugates of $H$ in $G$, we see easily that 
$G-\{1\}$ is partitioned into the sets $K-\{1\}, H^{g_1}-\{1\},\dots,
H^{g_m}-1$ and two elements picked from different sets in this 
partition will not commute. Thus it follows that the non-commuting complex 
based on the elements of $G-\{1\}$ decomposes as a join of the non-commuting 
complexes based on each set in the partition. To complete the picture 
one notices that each $H^g$ is isomorphic to $H$ and so contributes  
the same non-commuting complex as $H$ and furthermore since the conjugates of 
$H$ make up $G-K$, a simple count gives that 
$m=\frac{|G|-|K|}{|H|-1}=\frac{|H||K|-|K|}{|H|-1}=|K|$.

The sentence about $nc(G)$ follows from the fact that if we define 
$d(S)=\dime(S)+1$ for a simplicial complex $S$, then 
$d(S_1 * S_2)=d(S_1) + d(S_2)$. Thus since $d(NC(G))=nc(G)$ 
this proves the stated formula concering $nc(G)$.

Finally when $H$ and $K$ are abelian, $NC(H)$ and $NC(K)$ are just 
a set of points namely the nonidentity elements in each group. 
The short exact sequence of the join together with the fact that 
$NC(G)=BNC(G)$ is simply-connected finishes the proof.  
\end{proof}

\begin{ex}
$A_4$ is a Frobenius group with kernel $\mathbb{Z}/2\mathbb{Z} \times 
\mathbb{Z}/2\mathbb{Z}$ and complement $\mathbb{Z}/3\mathbb{Z}$. 
Thus proposition~\ref{pro: Frob} shows 
$NC(A_4) \simeq S^4 \vee S^4$ and $nc(A_4)=5$.
\end{ex}

\begin{claim}
$NC_2(A_5)$ is a $4$-spherical complex 
homotopy equivalent to a bouquet of $32$ $4$-spheres.
Thus it is $3$-connected and has odd Euler characteristic.
Hence $NC(A_5)$ has odd Euler characteristic.
\end{claim} 
\begin{proof}
The order of $A_5$ is $60=4 \cdot 3 \cdot 5$ of course. 
It is easy to check that there are five Sylow 2-subgroups $P$ 
which are elementary abelian 
of rank 2 and self-centralizing, i.e., $C_G(P)=P$, and are ``disjoint'', 
i.e., any two Sylow subgroups intersect only at the identity element.

Thus the picture for the vertices of $NC_2(A_5)$ is as 5 sets 
$\{S_i\}_{i=1}^5$ of 
size 3. (Since each Sylow 2-group gives 3 involutions.) 
Now since the Sylow 2-groups are self-centralizing, this means 
that two involutions in two different Sylow 2-subgroups, do not commute 
and thus are joined by an edge in $NC_2(A_5)$. Thus we see easily 
that $NC_2(A_5)$ is the join $S_1 * S_2 * S_3 * S_4 * S_5$.  

Using the short exact sequence for the join 
(see page 373, Exercise 3 in \cite{M}), 
one calculates easily that $NC_2(A_5)$ has the homology of a bouquet 
of $32$ $4$-spheres. Since the join of two path connected spaces 
$S_1 * S_2$ and $S_3 * S_4 * S_5$ is simply-connected, it follows 
that $NC_2(A_5)$ is homotopy equivalent to a bouquet of $32$ $4$-spheres, 
and since it is obviously $4$-dimensional, this completes all but 
the last sentence of the claim. The final sentence follows from 
proposition~\ref{pro: 2local} which says that $NC(A_5)$ has 
the same Euler characteristic as $NC_2(A_5)$ mod 2. 
\end{proof}

At this stage, we would like to make a conjecture:

\begin{conj}
\label{conj: contract}
If $G$ is a nonabelian simple group, then $NC(G)=BNC(G)$ 
has odd Euler characteristic.
\end{conj}

Recall the following famous theorem of Feit and Thompson:

\begin{thm}[Odd order theorem]
Every group of odd order is solvable.
\end{thm}

Notice, that if the conjecture is true, it would imply the odd order theorem.
This is because it is easy to see that a minimal counterexample 
$G$ to the odd order theorem would have to be an odd order nonabelian simple 
group. 
The conjecture would then say $BNC(G)$ has an odd Euler characteristic and 
proposition~\ref{pro: oddorder} would say $G$ was even order which 
is a contradiction to the original assumption that $G$ has odd order.

\section{General commuting structures}
Before we further analyze the $NC(G)$ complexes introduced in the last 
section, we need to extend our considerations to general commuting 
structures.

\begin{defn}
A commuting structure is a set $S$ together with a 
reflexive, symmetric relation $\sim$ on $S$. 
If $x, y \in S$ with $x \sim y$ we say $x$ and $y$ commute.
\end{defn}

\begin{defn}
Given a commuting structure $(S, \sim)$, the dual 
commuting structure $(S, \sim')$ is defined by \\  
\noindent
(a) $\sim'$ is reflexive and \\
\noindent
(b)For $x \neq y$, 
$x \sim' y$ if and only if $x$ does not commute with $y$ in 
$(S, \sim)$. \\ 

When it is understood 
we write a commuting structure as $S$ and its dual as $S'$.
It is easy to see that $S''=S$ in general.
\end{defn}

\begin{defn}
Given a commuting structure $S$, $C(S)$ is the 
simplicial complex whose vertices are the elements of $S$ 
and such that $[s_0,\dots,s_n]$ is a face of $C(S)$ if 
and only if $\{s_0,\dots,s_n\}$ is a commuting set i.e., 
$s_i \sim s_j$ for all $i$ and $j$. 
Similarly we define $NC(S)=C(S')$ and refer to the elements 
in a face of $NC(S)$ as a non-commuting set in $S$.
\end{defn}

Below are some examples of commuting structures which 
will be important in our considerations: \\
\noindent
(a) The nontrivial elements of a group $G$ form a commuting structure 
which we denote also by $G$, where $x \sim y$ if and only if $x$ and 
$y$ commute in the group $G$. In this case $C(G)$ and $NC(G)$ 
are the complexes considered in the previous sections. \\

\noindent
(b) The noncentral elements of a group $G$ form a commuting structure 
$G-Z(G)$ and $NC(G-Z(G))=BNC(G)$. Similarly the elements of order $p$ 
in $G$ form a commuting structure $G_p$ and $C(G_p)=C_p(G)$, 
$NC(G_p)=NC_p(G)$. \\

\noindent
(c) If $V$ is a vector space equipped with a bilinear map 
$[\cdot, \cdot]: V \otimes V \rightarrow V$. Then the nonzero elements of $V$ 
form a commuting structure also denoted by $V$, where $v_1 \sim v_2$ 
if $[v_1,v_2]=0$ or if $v_1 = v_2$. \\

\noindent
(d) In the situation in (c), we can also look at the set 
$P(V)$ of lines in $V$. The commuting structure on $V$ descends 
to give a well-defined commuting structure on $P(V)$, which 
we will call the projective commuting structure. \\

\noindent
(e) If $1 \rightarrow C \rightarrow G \rightarrow Q \rightarrow 1$ 
is a central extension of groups with $Q$ abelian then one can form 
$[\cdot, \cdot]: Q \rightarrow C$ by 
$$
[x,y]=\hat{x}\hat{y}\hat{x}^{-1}\hat{y}^{-1}
$$
where $x, y \in Q$ and $\hat{x}$ is a lift of $x$ in $G$ etc.
It is easy to see this is well-defined, independent of the choice of 
lift and furthermore that the bracket $[\cdot,\cdot]$ is bilinear. 
Thus by (c), we get a commuting structure on the nontrivial elements 
of $Q$ via this bracket. We denote this commuting structure by 
$(G;C)$. Notice in general it is not the same as the commuting 
structure of the group $G/C$ which is abelian in this case. 
More generally even if $Q$ is not abelian, one can define a 
commuting structure on the nontrivial elements of 
$Q$ from the extension by declaring 
$x \sim y$ if and only if $\hat{x}$ and $\hat{y}$ commute in $G$, 
we will denote this commuting structure by $(G;C)$ in general.

\begin{ex}
If $P$ is an extraspecial $p$-group of order $p^3$ then 
it has center $Z$ of order $p$ and $P/Z$ is an elementary abelian 
$p$-group of rank 2. Let $Symp$ denote the commuting structure 
obtained from a vector space of dimension two over $\mathbb{F}_p$ 
equipped with the symplectic alternating inner product $[x,y]=1$ 
where $\{x,y\}$ is a suitable basis. 
Then it is easy to check that $(P;Z)=Symp$. 
\end{ex}

One of the main results we will use in order to study the non-commuting 
complexes associated to commuting structures is a 
result of A. Bj\"orner, M. Wachs and V. Welker on ``blowup'' complexes which 
we describe next:

Let $S$ be a finite simplicial complex with vertex set $[n]$ (This 
means the vertices have been labelled 1,\dots,n). 
To each vertex $1\leq i \leq n$, we assign a positive integer $m_i$. 
Let $\bar{m}=(m_1,\dots,m_n)$, 
then the ``blowup'' complex $S_{\bar{m}}$ is defined as follows:

The vertices of $S_{\bar{m}}$ are of the form $(i,j)$ where 
$1\leq i \leq n, 1\leq j \leq m_i$. (one should picture 
$m_i$ vertices in $S_{\bar{m}}$ over vertex $i$ in $S$.)

The faces of $S_{\bar{m}}$ are exactly of the form 
$[(i_0,j_0),\dots,(i_n,j_n)]$ where $[i_0,\dots,i_n]$ is a face 
in $S$. (In particular $i_l \neq i_k$ for $l \neq k$.)

The result of Bj\"orner, Wachs and Welker describes $S_{\bar{m}}$ up 
to homotopy equivalence, in terms of $S$ and its links. 
More precisely we have:

\begin{thm}[Bj\"orner, Wachs, Welker~\cite{Bj}]
\label{thm: blowup}
For any connected simplicial complex $S$ with vertex set $[n]$ and 
 given $n$-tuple of positive integers $\bar{m}=(m_1,\dots,m_n)$ we 
have:

$$
S_{\bar{m}} \simeq S \vee \bigvee_{F \in S} [Susp^{|F|}(Lk(F))]^{\gamma(F)}.
$$
Here the $\vee$ stands for wedge of spaces, $Susp^k$ for $k$-fold 
suspension, $Lk(F)$ for the link of the face $F$ in $S$, $|F|$ for 
the number of vertices in the face $F$ (which is one more than the 
dimension of $F$) and 
$\gamma(F) = \prod_{i \in F} (m_i-1)$. Thus in the decomposition 
above, $\gamma(F)$ copies of $Susp^{|F|}(Lk(F))$ appear wedged together. 
\end{thm}

Note that in~\cite{Bj}, the empty face is considered a face in any complex. 
In the above formulation we are not considering the empty face as a face 
and thus have separated out the $S$ term in the wedge decomposition.

For example, if we take a 1-simplex $[1,2]$ as our complex $S$ 
and use $\bar{m}=(2,2)$, then it is easy to see that 
$S_{\bar{m}}$ is a circle. On the other hand all links in $S$ 
are contractible except the link of the maximal face $[1,2]$ which 
is empty (a ``$(-1)$-sphere''). Thus everything in the right hand 
side of the formula is contractible except the 2-fold suspension 
of this $(-1)$-sphere which gives a 1-sphere or circle as expected.

Also note that whenever some $m_i=1$, the corresponding $\gamma(F)=0$ 
and so that term drops our of the wedge decomposition.
Thus if $m_1=\dots=m_n=1$, the decomposition gives us nothing as 
$S_{\bar{m}}=S$.

Now for some examples of where this theorem applies:

\begin{cor}
\label{cor: center}
If $G$ is a finite group and $Z(G)$ is its center, then 
$BNC(G)$ is the blowup of $(G;Z(G))$ where each $m_i=|Z(G)|$. 
This is because everything in the same coset of the $Z(G)$ in $G$ 
commutes with each other and whether or not two elements 
from different cosets commute is decided in $(G;Z(G))$.
Thus we conclude 
$$
BNC(G) \simeq NC(G;Z(G)) \vee \bigvee_{F \in NC(G;Z(G))} 
[Susp^{|F|}(Lk(F))]^{(|Z(G)|-1)^{|F|}}.
$$
\end{cor}

Before we do say more, let us look at another example of  
theorem~\ref{thm: blowup} in our context. The proof is the same as 
that of example~\ref{cor: center} and is left to the reader.

\begin{cor}
\label{cor: vectorspace}
Let $(V,[\cdot,\cdot])$ be a vector space over a finite field $\mathbf{k}$ 
equipped with a nondegenerate bilinear form.
Then $NC(V)$ is a blowup of $NC(P(V))$ and we have 
$$
NC(V) \simeq NC(P(V)) \vee \bigvee_{F \in NC(P(V))}
[Susp^{|F|}(Lk(F))]^{(|\mathbf{k}|-2)^{|F|}}.
$$ 
\end{cor}

For example, it is easy to see that in $Symp$, if $[x,y]=0$ then 
$x$ is a scalar multiple of $y$. Thus in $P(Symp)$ two distinct 
elements do not commute and hence $NC(P(Symp))$ is a simplex. 
Thus all of its links are contractible except the link of the top face 
which is empty. This top face has $(p^2-1)/(p-1)=p+1$ vertices. 
Thus following example~\ref{cor: vectorspace}, we see 
$$
NC(Symp) \simeq \vee_{(p-2)^{(p+1)}} S^p.
$$

\begin{rem}
If $(V,[\cdot,\cdot])$ is a vector space over a finite field 
equipped with a nondegenerate symmetric inner product corresponding 
to a quadratic form $Q$, then one can restrict the commuting structure 
induced on $P(V)$ to the subset $S$ 
consisting of singular points of $Q$, i.e., lines $<x>$ with 
$Q(x)=0$. ($S$ will be nonempty only if $Q$ is of hyperbolic type.)
Among the non-commuting sets in $S$ are the 
subsets called ovoids defined by the property that every maximal singular 
subspace of $V$ contains exactly one element of the ovoid. (See \cite{G})
Thus these ovoids are a special subcollection of the facets of 
$NC(S)$.
\end{rem}

In general, we will find the following notion useful:

\begin{defn} 
If $(S,\sim)$ is a commuting set 
and $x \in S$, we define the centralizer of $x$ 
to be
$$
C(x)=\{ y \in S| y \sim x \}.
$$
\end{defn}

This allows us to define

\begin{defn}
If $(S, \sim)$ is a commuting set, we define an equivalence relation 
on the elements of $S$ by
$x \approx y$ if $C(x)=C(y)$. The equivalence classes are called 
the centralizer classes of $(S,\sim)$. 
We then define the core of $S$ to be the commuting set $(\bar{S}, \sim)$
where the elements are the centralizer classes of $S$ and the classes 
$[x]$ and $[y]$ commute in $\bar{S}$ 
if and only if the representatives $x$, $y$ 
commute in $S$.
\end{defn}

It is easy to see that $NC(S)$ is the blowup of $NC(\bar{S})$, 
where for the vertex $[x] \in NC(\bar{S})$ there are $m_x$ vertices 
above it in $NC(S)$, where $m_x$ is the size of the centralizer class $[x]$.
Thus once again theorem~\ref{thm: blowup} gives us:

\begin{thm}
\label{thm: core}
Let $S$ be a commuting set and $\bar{S}$ be its core, and suppose 
$NC(\bar{S})$ is connected. Then 
$$
NC(S) \simeq NC(\bar{S}) \vee \bigvee_{F \in NC(\bar{S})} 
[Susp^{|F|}(Lk F)]^{\gamma(F)}
$$
where $\gamma(F)=\prod_{[x] \subseteq F}(m_x - 1)$.
Here again $m_x$ is the size of the centralizer class $[x]$.
\end{thm}

\begin{rem}
Note that if $NC(S)$ is connected, $NC(\bar{S})$ will automatically 
be connected as it is the image of $NC(S)$ under a continuous map.
\end{rem}

With some abuse of notation, we will call $NC(\bar{S})$ the core of $NC(S)$.

Following the notation above, in a finite group $G$, the equivalence 
relation ``has the same centralizer group'' partitions $G$ into 
centralizer classes. The central elements form one class and the 
noncentral elements thus inherit a partition.

Notice if $[x]$ is the centralizer class containing $x$, then 
any other generator of the cyclic group $<x>$ is in the same class. 
Thus $[x]$ has at least $\phi(n)$ elements where $n$ is the order of 
$x$, and $\phi$ is Euler's totient function. Thus if $n>2$, then 
$[x]$ contains at least two elements. Also notice that every thing 
in the coset $xZ(G)$ is also in $[x]$ so if $Z(G) \neq 1$ we can 
also conclude $[x]$ contains at least two elements.

\begin{defn} A non-commuting set $S$ in $G$ is a nonempty 
subset $S$, such that the elements of $S$ pairwise do not commute.
A maximal non-commuting set $S$ is a non-commuting set which is not 
properly contained in any other non-commuting set of $G$.
\end{defn}

In general, not all maximal non-commuting sets of a group $G$ have 
the same size.

One obtains the following immediate corollary of theorem~\ref{thm: core}
(Assume $|G|>2$ for the following results):

\begin{cor}
\label{cor: dimension}
Let $G$ be a finite nonabelian group, and let $S_s$ 
denote the set of maximal non-commuting sets in 
$G$ of size $s$. Then, for $s >1$,
$$
rk(H_{s-1}(BNC(G))) \geq \sum_{F \in S_s} [\prod_{x \in F}
(1-\frac{1}{m_x})]. $$
where $m_x$ is the size of the centralizer class containing $x$.

In particular, if $G$ is an odd order group 
or if $G$ has a nontrivial center, then 
$$ {\tilde H}_{s-1}(NC(G)) \neq 0 $$ whenever $G$ has a maximal 
non-commuting set of size $s$.
\end{cor}
\begin{proof}
Let $X$ be the non-commuting complex associated to the core of $BNC(G)$.
Let us define a facet to be a face of a simplicial complex which is 
not contained in any bigger faces. Thus the facets of $BNC(G)$ consist 
exactly of maximal non-commuting sets of noncentral elements in $G$.

The first thing to notice is 
to every facet $F$ of $BNC(G)$, there corresponds a facet $\bar{F}$ 
of $X$, and further more 
this correspondence preserves the dimension of the facet (or equivalently 
the number of vertices in the facet).

Since the link of a facet is always empty (a $(-1)$-sphere),
in the wedge decomposition of theorem~\ref{thm: core}, 
we get the suspension of this empty link as a contribution. 
If the facet $F$ has $n$ vertices in it, then we suspend $n$ times to get 
a $(n-1)$-sphere. Thus to every maximal non-commuting set of size $s$, 
we get a $s-1$ sphere contribution from the corresponding facet in $X$.
In fact we get $\gamma(F)$-many such spheres from the facet $F$.
However above each facet $\bar{F}$ in $X$, there correspond 
$\prod_{[x] \in \bar{F}}(m_x)$ many facets in $BNC(G)$. 
Thus in the sum over the facets of $BNC(G)$ stated in the theorem, 
we divide $\gamma(F)$ by $\prod_{x \in F}(m_x)$ in order to count 
the contribution from the facet $\bar{F}$ in $X$ the correct number of times.

Notice that $\frac{\gamma(F)}
{\prod_{x \in F}(m_x)}=\prod_{x \in F}(1-\frac{1}{m_x}) 
\geq \frac{1}{2^{|F|}}$ 
if all the centralizer classes have size bigger than one.
So if $G$ has the property that the size of the 
centralizer classes of noncentral 
elements is always strictly bigger than one, for example if $G$ is odd 
or if $G$ has a nontrivial center, then 
$$
2^s rk(H_{s-1}(BNC(G))) \geq |S_s|
$$
for all $s \in \mathbb{N}$. 
In particular $H_{s-1}(BNC(G)) \neq 0$ whenever $G$ has a maximal 
non-commuting set of size $s$. Also observe that 
${\tilde H}_0 (NC(G)) \neq 0$ whenever $G$ has maximal non-commuting set
of size $1$, i.e. a singleton consisting of a nontrivial central element
(except for the trivial case when $G$ has order 2.)
\end{proof} 

It is easy to see from this proof that a $p$-local
version of corollary \ref{cor: dimension} is also true. 
We state here only the last part of this result for odd primes:
\begin{cor}
\label{cor: dimension2}
Let $p$ be an odd prime. Then,
$$ {\tilde H}_{s-1}(NC_p(G)) \neq 0 $$ whenever $G$ has a maximal
non-commuting $p$-set of size $s$.
\end{cor}

The same is true for $p=2$ under the additional condition $2||Z(G)|$.

\begin{rem} Notice that the conclusion of corollary \ref{cor:
dimension} is
consistent with the simple connectedness of $BNC(G)$, because there is no
maximal non-commuting set of size $2$. To see this, observe that
whenever there is a non-commuting  set $\{ a, b\}$ with two elements, we
can form a bigger non-commuting set $\{ a, b, ab \}$. 

Also notice that this is no
longer the case for the $p$-local case. For example, when $G=D_8 = 
\langle a, b | a^2= b^2=c^2=1, [a,b]=c, \ c \ central \rangle$, the complex 
$BNC_2(G)$ is a rectangle with vertices $a,b,ac,bc$ which is
the inflated complex corresponding to maximal $2$-set $\{ a, b \}$.
In particular, $BNC_2(G)$ is not simply connected in general. 
\end{rem}

\begin{rem}
One of the things that corollary~\ref{cor: dimension} says is that if one 
wants to calculate $nc(G)$, 
the answer which is obviously the dimension of $BNC(G)$ plus 
one can also be determined by finding the highest nonvanishing homology 
of $BNC(G)$ in the case when $Z(G) \neq 1$ or $G$ is odd order. 
Thus the answer is determined already by the homotopy 
type of $BNC(G)$ in this situation. If $Z(G)=1$ and $G$ has even order, 
this is no longer 
true, for example $NC(\Sigma_3)$ is contractible and so does not 
have any positive dimensional homology.
\end{rem}

Sometimes one can show the non-commuting complex for the 
core of $(S,\sim)$ is contractible as 
the next lemma shows:

\begin{lem}
\label{lem: coregoes}
If $(S,\sim)$ has a centralizer class $[x]$ where $[x]=C(x)$ 
then if the core is $\bar{S}$, then $NC(\bar{S})$ is contractible.
\end{lem}
\begin{proof}
This is because it is easily seen that $NC(\bar{S})$ is a cone on 
the vertex $[x]$ as everything outside $[x]$ does not commute 
with $x$ as $[x]=C(x)$.
\end{proof} 

Thus for example:

\begin{ex}
Let $p$ be a prime, then the core of $BNC(\Sigma_p)$ is contractible.
\end{ex}
\begin{proof}
The cycle $x=(1,2,\dots,p)$ in $\Sigma_p$ has $C(x)=<x>$ by a 
simple calculation. Thus $C(x)=[x]$ and so the result follows 
from lemma~\ref{lem: coregoes}.
\end{proof} 

We now study an important special case:

\begin{defn} 
A $TC$-group $G$ is a group where the commuting relation on the 
noncentral elements is transitive. This is equivalent to the condition 
that all proper centralizer subgroups $C(x) \subset G$ are abelian.

For example any minimal nonabelian group like $S_3$, $A_4$ or 
an extraspecial group of order $p^3$.
\end{defn}

\begin{cor}
\label{cor: equiv}
If $(S, \sim)$ is a commuting set where $\sim$ is also transitive (i.e., 
$\sim$ is an equivalence relation), then 
$NC(S)$ is homotopy equivalent to a wedge of spheres of the same 
dimension. (We allow the ``empty'' wedge, i.e., we allow the case $NC(S)$ 
is homotopy equivalent to a point). The dimension of the spheres 
is equal to $n-1$ where $n$ is the number of equivalence classes in $(S,\sim)$ 
and the number of spheres appearing is $\prod_{i=1}^n (m_i-1)$ where 
$m_i$ is the size of equivalence class $i$.

Thus if $G$ is a $TC$-group, then 
$BNC(G)$ is homotopy equivalent to a wedge of spheres of dimension 
$nc(G)-1$ and the number of spheres is given by a product as above 
where the $m_i$ are the orders of the distinct proper centralizer 
groups of $G$.
\end{cor}
\begin{proof}
If $\sim$ is an equivalence relation, it is easy to see that 
the centralizer classes are exactly the $\sim$ equivalence classes. 
Thus one sees that the non-commuting complex for the core, 
$NC(\bar{S})$, is a simplex. Thus in Theorem~\ref{thm: core}, 
all terms drop out except that corresponding to the maximum face 
in $NC(\bar{S})$  
where the link is empty. This link is suspended to give a sphere 
of dimension equal to the number of equivalence classes minus one. 
The number of these spheres appearing in the wedge decomposition is 
the product $\prod_{i=1}^n(m_i-1)$ where $m_i$ is the size of 
equivalence class $i$ and the product is over all equivalence 
classes. (Thus this  
can be zero if one of the equivalence classes has size one, in 
which case the complex is contractible).  

In the case of $BNC(G)$, for $G$ a $TC$-group, one just has to 
note that $C(x)$ is the centralizer class $[x]$ for any 
noncentral element $x$. 
\end{proof}

\begin{rem}
In Corollary~\ref{cor: equiv}, one did not actually have to use 
the general result of Bjorner, Wachs and Welker since it is easy to see 
that in this situation, $NC(S)$ is the join of each equivalence class 
as discrete sets and a simple count gives the result.
\end{rem}

Now notice in the case that $(S,\sim)$ has $\sim$ transitive, 
the above analysis shows that $NC(S)$ is a join of discrete sets. 
Thus it is easy to see that $NC(S)$ is shellable (This is because the 
facets of $NC(S)$ are just sets where we have chosen exactly 
one element from each of the $\sim$ equivalence classes. We can linearly order 
the equivalence classes and then lexicographically order the facets. It 
is easy to check that this is indeed a shelling.)

Given a shelling of a simplicial complex, there are many combinatorial 
equalities and inequalities which follow (See~\cite{S}.)
Since these are not so deep in the above general context, 
we will point out only the interpretation when applied to $BNC(G)$. 
Recall pure shellable just means shellable where all the 
facets have the same dimension.

\begin{pro}
\label{pro: Shell}
If $G$ is a nonabelian group such that $BNC(G)$ is pure shellable, 
e.g., $G$ a $TC$-group like $\Sigma_3$ or $A_4$, 
then if one sets $nc_i$ to be the number of non-commutative sets of 
non-central elements which have size $i$, one has:
$$
(-1)^jC(nc(G),j) + \sum_{k=1}^j (-1)^{j-k} C(nc(G)-k,j-k)nc_k \geq 0
$$
for all $1 \leq j \leq nc(G)$ where $C(n,k)$ is the usual 
binomial coefficient.
\end{pro}
\begin{proof}
Follows from a direct interpretation of the inequalities in 
\cite{S}, page 4, Theorem 2.9. One warning about the notation in 
that paper is that $|\sigma|$ means the number of vertices in $\sigma$ 
and the empty face is considered a simplex in any complex.
\end{proof}

Using that $nc_1=|G|-|Z(G)|$, $nc_2=\frac{|G|}{2}(|G|-m)$ 
where $m$ is the number of conjugacy classes in $G$, one gets for example 
from the inequality with $j=2$ above:

\begin{cor}
\label{cor: Shell}
Let $G$ be a nonabelian group with a transitive commuting relation,
i.e., if $[g,h]=[h,k]=1$, then $[g,h]=1$ for every noncentral $g,h,k \in G$.
Then, 
$$
nc(G)(nc(G)-1) + |G|(|G|-m) - 2(nc(G)-1)(|G|-|Z(G)|) \geq 0.
$$
where $m$ denotes the number of conjugacy classes in $G$.
\end{cor}

\section{Duality}
Let $(S, \sim)$ be a finite set with a commuting relation. Suppose  
the commuting complex $C(S)$ breaks up as a disjoint union of 
path components $C(S_1), \dots, C(S_n)$ where of course we are using 
$S_i$ to stand for the vertex set of component $i$. 

Then notice in the corresponding non-commuting complex, $NC(S)$ 
we have 
$NC(S)=NC(S_1) * \dots * NC(S_n)$ where $*$ stands for the join operation 
as usual.

We state this simple but useful observation as the next lemma:

\begin{lem}[Duality]
\label{lem: duality}
Let $(S, \sim)$ be a commuting set then if 
$$
C(S) = \bigsqcup_{i=1}^n C(S_i)
$$
where $\bigsqcup$ stands for disjoint union, then we have
$$
NC(S) = *_{i=1}^n NC(S_i)
$$
where $*$ stands for join.
\end{lem}
 
Thus in some sense ``the less connected $C(S)$ is, the more connected 
$NC(S)$ is.''

We can apply this simple observation to say something about 
the complexes $NC_p(G)$ in general. 

\begin{thm}
\label{thm: Quillen}
Let $G$ be a finite group and $p$ a prime such that $p | |G|$. 
Let $P$ be a Sylow $p$-group of $G$ and define $N$ 
to be the subgroup generated by $N_G(H)$ as $H$ runs over all 
the nontrivial subgroups of $P$. 

Then $NC_p(G)$ is $(|G:N|-2)$-connected and in fact it is 
the $|G:N|$-fold join of some complex with itself.
\end{thm}
\begin{proof}
By Quillen~\cite{Q}, if $S_1, \dots, S_n$ are the components 
of $A_p(G)$, then under the $G$-action, $G$ acts transitively 
on the components with isotropy group $N$ under suitable choice 
of labelling. Thus the components are all simplicially equivalent 
and there are $|G:N|$ many of them. 

However we have seen that $A_p(G)$ is $G$-homotopy equivalent to 
$C_p(G)$ and so we have the same picture for that complex. 
Thus $C_p(G)$ is the disjoint union of $|G:N|$ copies of some 
simplicial complex $S$. Thus by lemma~\ref{lem: duality}, 
$NC_p(G)$ is the $|G:N|$-fold join of the dual of $S$ with itself. 
To finish the proof one just has to note that a $k$-fold join of  
nonempty spaces is 
always $(k-2)$-connected.  
\end{proof}

\bigskip

\noindent
Dept. of Mathematics \\
University of Wisconsin-Madison, \\
Madison, WI 53706, U.S.A. \\
E-mail address: pakianat@math.wisc.edu \\

\bigskip 

\noindent
Dept. of Mathematics \& Statistics \\
McMaster University \\
Hamilton, ON, Canada \\
L8S 4K1. \\
E-mail address: yalcine@math.mcmaster.ca \\

\end{document}